# UN PROBLEMA DA DISCUTERE
## Una rappresentazione geometrica del teorema del coseno

*Dette a, b, c le lunghezze dei lati di un qualunque triangolo e detta γ l'ampiezza dell'angolo opposto al lato c, vale l'uguaglianza $a^2 + b^2 = c^2 + 2ab \cos \gamma$.*

L'enunciato è il noto *teorema del coseno*. Talvolta si parla anche di teorema di Carnot, dal nome di Lazare Carnot (1753-1823), ma in realtà il teorema era noto ben prima di Carnot, almeno dal XVI secolo (qualcuno lo attribuisce a François Viète). In termini geometrici l'enunciato si ritrova addirittura in Euclide, ovviamente senza un ricorso esplicito a funzioni trigonometriche. In effetti, la Proposizione 13 del Libro 2 degli *Elementi* afferma quanto segue (la Proposizione 12 è analoga, ma si riferisce ai triangoli ottusangoli).

**Proposizione 13.** *Nei triangoli acutangoli il quadrato sul lato opposto all'angolo acuto è minore dei quadrati sui lati che comprendono l'angolo acuto per il doppio del rettangolo compreso da uno dei lati intorno all'angolo acuto, quello su cui cade la perpendicolare, e dalla retta staccata dalla perpendicolare all'angolo acuto.*

Il testo ci appare assolutamente oscuro nella parte finale; ma c'è una spiegazione con riferimento ad una figura. Dicendo «*il quadrato […] è minore dei quadrati […] per il doppio del rettangolo*» si intende dire che il doppio del rettangolo è la differenza fra la somma dei quadrati costruiti su due lati (quelli che comprendono l'angolo acuto) e il quadrato costruito sul terzo lato. Euclide rappresenta geometricamente tale differenza. Se *AC* (figura 1) è il lato opposto all'angolo acuto in *B*, il rettangolo che rappresenta la metà della differenza citata ha per dimensioni *BC* e *BD* (*BC* è «*uno dei lati intorno all'angolo acuto*», mentre *BD* è il segmento «*staccato dalla perpendicolare all'angolo acuto*»).

Confondendo un segmento con la corrispondente misura, si ha quindi

$$AC^2 = AB^2 + BC^2 - 2\,BC \cdot BD.$$

Per ritrovare l'enunciato usuale è sufficiente osservare che $BD = AB \cos \beta$.

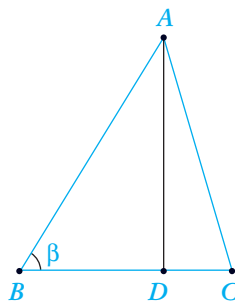

Figura **1**

Sono note varie dimostrazioni del teorema; qui propongo una recente dimostrazione di Al Cuoco, americano di origini italiane, autorevole studioso di didattica della matematica e direttore del *Center for Mathematics Education* presso l'*Education Development Center* nel Massachusetts (USA); la dimostrazione è riportata in [H].

Supponiamo, in un primo tempo, che il triangolo sia acutangolo.
Costruiamo i quadrati sui tre lati del triangolo, esternamente al triangolo; quindi tracciamo le tre altezze del triangolo e prolunghiamole in modo che ogni quadrato sia diviso in due rettangoli (fig. 2).
Calcoliamo ora le aree dei due rettangoli $R_1$ ed $R_2$, colorati in figura 1, che hanno per dimensioni un lato del triangolo e la proiezione di un secondo lato sul primo.

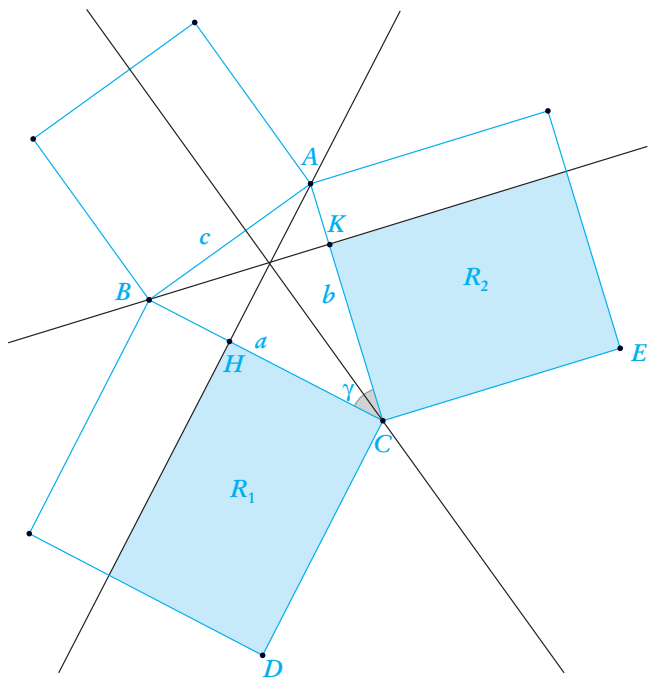

Figura **2**

Se si conoscono le nozioni base della trigonometria il calcolo delle aree è semplice:

$$\text{Area } R_1 = CD \times CH = a \times b \cos \gamma$$
$$\text{Area } R_2 = CE \times CK = b \times a \cos \gamma$$

Si conclude così che i due rettangoli sono equivalenti.

Allo stesso risultato si arriva rapidamente anche senza trigonometria: basta riconoscere la similitudine dei triangoli *CAH* e *CBK* (triangoli rettangoli con l'angolo in *C* in comune). Si ha allora la proporzione

$$CA:CH = CB:CK \qquad \text{da cui} \qquad b \times CK = a \times CH.$$

In modo analogo si trovano altre due coppie di rettangoli equivalenti, come suggerito dalla figura 3.

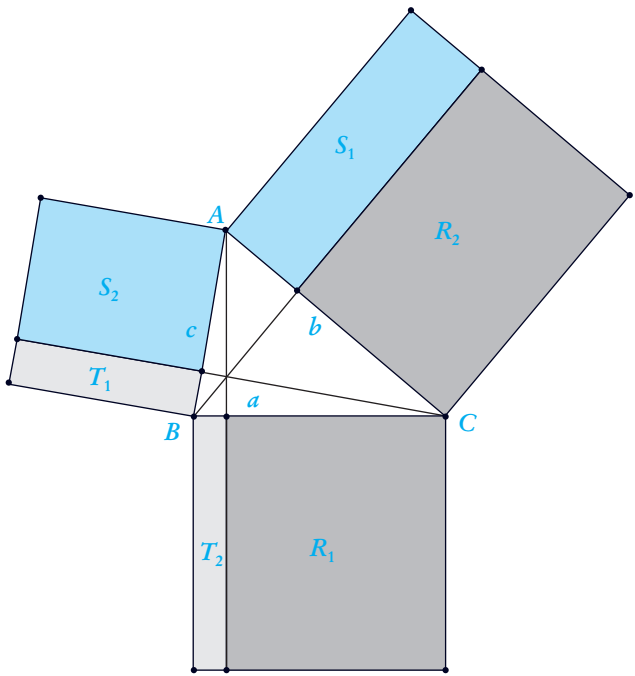

Figura **3**

Abbiamo così scomposto i tre quadrati in sei rettangoli, a due a due equivalenti. Il ragionamento è ora facile:

$$a^2 = R_1 + T_2 = R_2 + T_1 = (b^2 - S_1) + (c^2 - S_2) = b^2 + c^2 - 2S_1$$

e, ricordando che $S_1 = bc \cos \alpha$, si conclude appunto

$$a^2 = b^2 + c^2 - 2bc \cos \alpha.$$

Se il triangolo è *rettangolo* le altezze coincidono con i cateti e si ottiene la classica configurazione del I teorema di Euclide, da cui segue immediatamente il teorema di Pitagora.

Se il triangolo è ottusangolo, il ragionamento funziona ancora, a patto di accettare aree *negative*. Con riferimento alla figura 4, il quadrato costruito sul lato *AC* geometricamente non è la somma, ma è la *differenza* dei due rettangoli $R_2$ = *CEFI* ed $S_1$: in effetti, $R_2$ è tutto il rettangolo *CEFI* e, togliendo $S_1$ da $R_2$, rimane appunto il quadrato costruito sul lato. Un discorso analogo vale per il quadrato costruito sul lato *AB*. Come nel caso precedente, $R_2$ è equivalente ad $R_1$ ed $S_1$ è equivalente ad $S_2$.

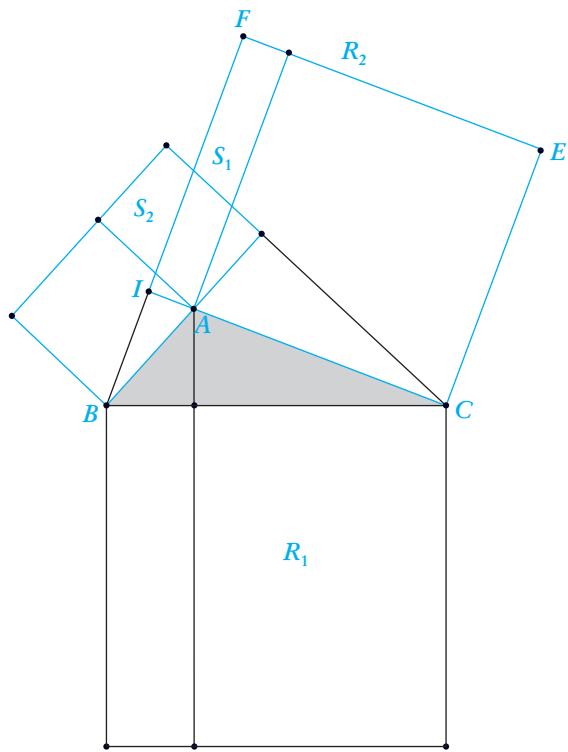

Figura **4**

Ancora un'osservazione.
Dati tre numeri *L*, *M*, *N*, il sistema

$$\begin{cases} x + y = L \\ x + z = M \\ y + z = N \end{cases}$$

ammette sempre una e una sola soluzione.

Se *L*, *M*, *N* sono le aree dei tre quadrati costruiti sui lati di un triangolo, la soluzione è rappresentata geometricamente dalle aree dei rettangoli considerati in figura 3.

Se *L*, *M*, *N* sono invece le lunghezze dei lati di un triangolo, allora *x*, *y*, *z* sono le lunghezze dei segmenti in cui i lati stessi sono divisi dai punti di tangenza con il cerchio inscritto (figura 5). In tal caso, siccome ciascuno dei tre numeri *L*, *M*, *N* è minore della somma degli altri due, si ottengono come soluzione tre numeri positivi (infatti, la soluzione è data dalla terna $(L+M-N)/2$, $(L+N-M)/2$, $(M+N-L)/2$). Nel caso dei quadrati la soluzione è formata da tre numeri positivi se e solo se il triangolo è acutangolo.

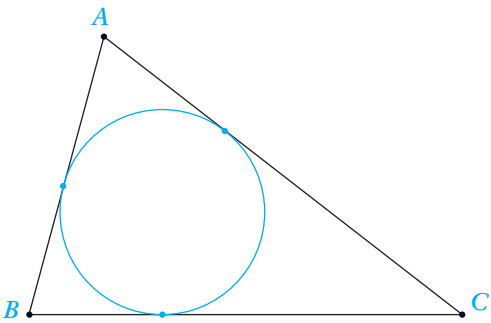

Figura **5**

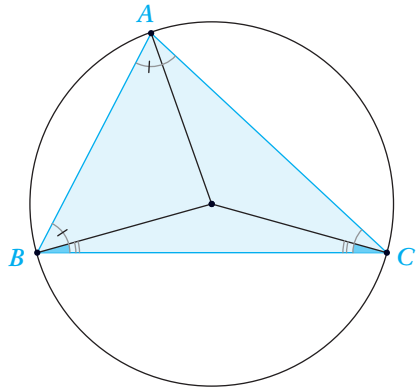

Figura **6**

Un'interpretazione geometrica è facile anche quando *L*, *M*, *N* sono le ampiezze degli angoli di un triangolo. Infatti, ogni angolo del triangolo è diviso in due parti da un raggio del cerchio circoscritto e la soluzione del sistema corrisponde proprio alle ampiezze di queste parti (figura 6), perché i tre triangoli che si ottengono congiungendo i vertici con il centro del cerchio sono ovviamente isosceli. Si osservi che anche in quest'ultimo caso la soluzione è formata da tre numeri positivi se e solo se il triangolo è acutangolo.

**Claudio Bernardi**
claudio.bernardi@uniroma1.it

**Bibliografia**

[H] R.E. Howe, *The Cuoco Configuration*, «The American Mathematical Monthly», Vol. 120 (2013), pagg. 916-923.